\documentclass[10pt]{pja00}

\usepackage{eucal}
\usepackage{amscd}
\usepackage[all,cmtip]{xy}
\usepackage{rotating}
\usepackage{hyperref,mathrsfs}
\usepackage{tikz}

\usepackage{amsmath,amssymb}

\newtheorem{theorem}{Theorem }[section]

\newtheorem{observation}[theorem]{Observation}

\def\I{\mathbf{I}}

\def\PG{\mathbf{PG}}

\newcommand{\Aut}{\mathrm{Aut}}

\newcommand{\mE}{\mathcal{E}}

\newcommand{\mP}{\mathcal{P}}

\newcommand{\mM}{\mathcal{M}}

\newcommand{\hW}{\mathcal{W}}
\newcommand{\hH}{\mathcal{H}}

\newcommand{\mB}{\mathcal{B}}

\runninghead{Classification of STGQs, I}
\title{Classification of STGQs, I}

\Subject[2000]{}

\Author{1}{Koen}{Thas}

\affiliation{1}{{Ghent University},
{Department of Mathematics},
{Krijgslaan 281, S25, B-9000 Ghent, Belgium,}
\texttt{koen.thas@gmail.com}}

\date{}

\begin{document}
\maketitle

\begin{abstract}
We describe new classification results in the theory of generalized quadrangles (= Tits-buildings of rank $2$ and type $\mathbb{B}_2$), more precisely in the (large) sub theory of  skew translation generalized quadrangles (``STGQs''). Some of these involve, and solve, long-standing open problems. 
\end{abstract}

%\setcounter{tocdepth}{1}
%\tableofcontents

\bigskip
\section{Introduction}

Generalized $n$-gons were introduced by Tits in a famous work on triality \cite{Ti} of 1959,  in order to propose an axiomatic and combinatorial treatment for  semisimple algebraic groups (including Chevalley groups  and groups of Lie type) of relative rank $2$. They are the central rank $2$ incidence geometries, and the atoms of the more general ``Tits-buildings.'' 
If the number of elements of a generalized $n$-gon is finite, a celebrated result of Feit and Higman \cite{FeHi} guarantees that $n$ is restricted to the set $\{ 3,4,6,8\}$.

Note that projective planes are nothing else than generalized $3$-gons. Generalized $4$-gons are also called {\em generalized quadrangles}, and certainly in the finite case, they are considered as being the main players in the class of generalized $n$-gons, and one of the most studied types of incidence geometries.

The most fruitful way to construct finite generalized quadrangles is through a now standard group coset geometry construction, in which 
a group $E$ provided with certain sets of subgroups $\mE = \{ E_i \vert i \in I \}$  $\mE^* = \{ E_i^* \vert i \in I \}$ and  satisfying some strong intersection properties, is used to represent a generalized quadrangle. Such a system of groups $(\mE,\mE^*)$ is called a {\em Kantor family} for $E$, and the defining properties are as follows.

\begin{itemize}
\item
For some $s, t \in \mathbb{N} \setminus \{0,1\}$, $\vert I \vert = t + 1$, $\vert E \vert = s^2t$, each $E_i$ has order $s$ and each $E_j^*$ has order $st$.
\item
For each $i$ (in $I$), we have $E_i \leq E_i^*$.
\item
For distinct $i, j$ and $k$ (in $I$), we have $E_iE_j \cap E_k = \{ \mathrm{id}\}$.
\item
For distinct $i$ and $j$ (in $I$), we have $E_i^* \cap E_j = \{ \mathrm{id}\}$.
\end{itemize}

Generic points of the quadrangle
are elements of the group, generic lines are left (or right | the choice does not change the isomorphism class of the geometry) cosets of the subgroups of type $E_i$, and some other special points and lines exist.
A fundamental feature of this construction is that $E$ naturally acts as an automorphism group of the geometry, sharply transitively on the generic points. See \cite{LEGQ} for the details. Call generalized quadrangles with such a group coset representation ``elation generalized quadrangles'' (``EGQs''); the group $E$ is the {\em elation group},
and there is a special point $(\infty)$ through which all the lines are fixed by $E$ (the ``elation point'').

In the long literature of generalized quadrangles, it has appeared that one special type of EGQ plays a central role. Such EGQs are called {\em skew translation generalized quadrangles} (``STGQs''), and they are specializations of EGQs in the sense that $(\infty)$ satisfies an additional combinatorial property called ``regularity." Except for the quadrangles associated to Hermitian varieties in $4$-dimensional projective space, all known classes of generalized quadrangles are STGQs, up to a combination of duality and Payne-integration. This observation strongly motivates the necessity to understand STGQs, and perhaps aiming for classification results in every which way. Only very partially, mostly in an influential paper by Payne \cite{STGQ89}, such results have been obtained for particular types of STGQs up till recently | see also \S \ref{HPMGQ}.

This classifcation program is the main purpose of this note.

In fact, more precisely, we want to understand the category $\mathbf{S}$, where objects are triples $(\Gamma,x,E)$, with $\Gamma$ an STGQ, and $x$
a regular elation point with respect to an elation group $E$. Morphisms are natural: If $A = (\Gamma,x,E)$ and $A' = (\Gamma',x',E')$, then 
$\mathrm{Hom}(A,A')$ consists of morhisms $\gamma: \Gamma \longrightarrow \Gamma'$ which map $x$ to $x'$ and $E$ to $E'$ (the 
latter meaning that for $e \in E$, $\gamma \circ e \in E'$). For now, we imagine the STGQs being finite (although the infinite case should eventually be handled as well).
Several subtleties arise:
\begin{itemize}
\item[(a)]
{\em one} fixed STGQ $\Gamma$ could have {\em several} (collinear) elation points $x$ and $x'$ (if the points would be non-collinear, a classification of such STGQs is known \cite{KT-HVM-TAMS}) - it can be easily shown however that there {\em always} exists automorphisms $\nu$ of $\Gamma$ mapping $x$ to $x'$;
\item[(b)]
even if one fixes the couple $(\Gamma,x)$, by recent work of Rostermundt and independently the author of this paper, see \cite{Rost,Basic}, it is known that examples exist which have {\em different}, even {\em non-isomorphic} elation groups $E$ and $E'$. As $E \not\cong E'$, $(\Gamma,x,E)$ and $(\Gamma,x,E')$ cannot be elements of the same isomorphism class in $\mathbf{S}$.
\end{itemize}

So classifying isomorphism classes in $\mathbf{S}$ is a strictly finer job than ``just'' classifying isomorphism classes of STGQs.

The examples of STGQs in (b) arise as the rational points and lines on a Hermitian variety in $\mathbf{PG}(3,q^2)$ ($q$ a power of $2$) (the point $x$ is arbitrary due to the transitivity of its automorphism group). This quadrangle is denoted by $\mathcal{H}(3,q^2)$.
Payne asked in 2004 (cf. \cite{StanPriv}) whether these examples are the {\em only (finite) STGQs} with different elation groups for the same point. We will come back to this question in the last section of this paper.

\bigskip
\section{Explanatory definitions}

In this paper,
a {\em generalized quadrangle} (``GQ'') is a point-line incidence geometry $\Gamma = (\mP,\mB,\I)$ for which  the following axioms are satisfied:

\begin{itemize}
\item[(i)]
$\Gamma$ contains no ordinary $k$-gon (as a subgeometry), for $2 \leq k < 4$;
\item[(ii)]
any two elements $x,y \in \mP \cup \mB$ are contained in some ordinary $4$-gon in
$\Gamma$;
\item[(iii)]
there exists an ordinary $5$-gon in $\Gamma$.
\end{itemize}

Here, $\mP$ is the point set, $\mB$ the line set, both non-empty and disjoint, and $\I$ is a symmetric relation on $(\mP \times \mB) \cup
(\mB \times \mP)$ called ``incidence,'' which tells us how points and lines are related. So $y \I Y$, with $y \in \mP$ and $Y \in \mB$, means that 
$y$ is incident with $Y$ (and $Y$ is incident with $y$).

An ``ordinary $3$-gon'' (e.g.) is a set of three points, two by two collinear but not all incident with the same line.

By (iii), generalized quadrangles have at least three points per line and three lines per point.

\subsection{Duality}
Note that points and lines play the same role in the axioms; this is the principle of ``duality.''

\subsection{Order}
All generalized quadrangles have an {\em order} $(s,t)$; there exist constants $s, t$ such that 
the number of points incident with a line is $s + 1$, and the number of lines incident  with a point is 
$t + 1$, cf. \cite{PT}. See also \cite{Order} for a detailed discussion regarding parameters of generalized $n$-gons.

Note that an ordinary quadrangle is just a ``generalized quadrangle without (iii),'' of order $(1,1)$ | we call such a subgeometry also ``apartment'' (of $\Gamma$).

\subsection{Automorphisms}
An {\em automorphism} of a generalized quadrangle $\Gamma = (\mP,\mB,\I)$ is a bijection of $\mP \cup \mB$ which preserves $\mP$, $\mB$ and incidence.
The full set of automorphisms of a GQ forms a group in a natural way | the {\em automorphism group} of $\Gamma$, denoted $\Aut(\Gamma)$. It is one of its most important invariants.

\bigskip
\label{HPMGQ}
\section{The Moufang condition}

In \cite{Titslect}, Tits proved roughly that there is a one-to-one correspondence between buildings of irreducible spherical type and rank $r \geq 3$, and the algebraic absolutely simple groups of relative rank $r$. In order to have a similar statement in the rank $2$ case | the case of generalized $n$-gons that is to say | one must impose an extra condition, called the {\em Moufang condition}. We describe it for generalized quadrangles.

Let $A$ be an apartment of a GQ $\Gamma$. A {\em root} $\gamma$ of $A$ is a set of $5$ different elements $e_0,\ldots,e_4$ in $A$ such that $e_i \I e_{i + 1}$ (where the 
indices are taken in $\{0,1,2,3,4\}$), and $e_0, e_4$ are the {\em extremal elements} of $\gamma$. There are two types of roots, depending on whether the extremal elements are lines or points; in the second case we speak of {\em dual roots} to make a distinction between the types. Also, a (dual) root $\gamma$ without its extremal elements | the {\em interior} of $\gamma$ | is denoted by $\dot{\gamma}$ and called (dual) {\em i-root}.

If $\mM$ is a subgeometry of $\Gamma$, by $\Aut(\Gamma)^{[\mM]}$ we denote the subgroup of the automorphism group $\Aut(\Gamma)$ of $\Gamma$ which fixes every line incident with a point of $\mM$ and every point incident  with a line of $\mM$. Now a (dual) root $\gamma$ is {\em Moufang} if $\Aut(\Gamma)^{[\dot{\gamma}]}$ acts transitively on the apartments containing 
$\gamma$. In fact,  $\Aut(\Gamma)^{[\dot{\gamma}]} =: A(\dot{\gamma})$ then acts {\em sharply} transitively on these apartments.
Once a (dual) root $\gamma$ is Moufang, all  (dual) roots with interior $\dot{\gamma}$ are also Moufang, with respect to the same group $A(\dot{\gamma})$. (The latter groups are uniquely defined by $\dot{\gamma}$ and the Moufang property.)
In a natural way, we also
use the terms ``Moufang i-root'' and ``dual Moufang i-root'', and the elements of $A(\dot{\gamma})$ are called {\em root-elations}.

Now $\Gamma$ is {\em half Moufang} if all roots or all dual roots are Moufang. It is {\em Moufang} if all roots {\em and} dual roots are.

All Moufang generalized quadrangles were classified in the classical work \cite{TiWe}, as the main and hardest step in the classification of all Moufang buildings of rank $2$ | see \cite{TiWe} for a historical sketch. In the finite case, this already followed from work of Fong and Seitz \cite{BN1,BN2}, see also Chapter 8-9 of \cite{PT}. Historical details can be found in the survey \cite{KT-HVM}; see also \cite{JAT-KT-HVM}.

In \cite{PT} (and the references therein), the strength of the local automorphic theory for generalized quadrangles became clear: to eventually end up with a geometric treatment of the classification of Moufang quadrangles, the authors developed a theory which studies local Moufang conditions for quadrangles, and eventually the theory reached far beyond the eventual goal. One instance is the theory of EGQs | see Chapter 8 of \cite{PT} and 
the recent monograph \cite{LEGQ}. Other instances of local Moufang theory are surveyed in, e.g., \cite{KT-HVM}, see also the more recent paper \cite{KTlocal}.

\bigskip
\section{Skew translation quadrangles}

An important class of STGQs consists of those STGQs (called ``flock quadrangles,'' and introduced in \cite{flocks87}) which arise through an intricate construction as a group coset geometry, from a 
{\em flock} of a quadratic cone $\mathcal{K}$ in $\mathbf{PG}(3,q)$ ($\mathbb{F}_q$ a finite field). A {\em flock} is a partition of  the $\mathbb{F}_q$-rational points of $\mathcal{K} \setminus \{ \mbox{vertex}\}$ into $q$ disjoint irreducible conics.

In the fundamental paper \cite{STGQ89}, Payne studied local Moufang conditions in GQs, partly to understand and generalize   flock quadrangles. This motivated him to introduce {\em skew translation generalized quadrangles} (``STGQs'').  One can prove that if $(\Gamma,x,E)$ is an STGQ of finite 
order $(s,t)$,  all dual i-roots on $x$ are {\em Moufang for the same group} $\mathbb{S} \leq E$; this means that $\mathbb{S}$ has size $t$, and that 
all its elements fix every point collinear with $x$. The elements of $\mathbb{S}$ are called {\em symmetries} with center $x$. The existence of 
$\mathbb{S}$ forces $x$ to have a certain combinatorial property called {\em regularity}, and conversely, one can show that an EGQ with 
regular elation point is an STGQ.

In \cite{STGQ89}, Payne introduced and studied a particular class of STGQs, called ``MSTGQs.'' They are STGQs $(\Gamma,x,E)$ with the following properties.
\begin{itemize}
\item[(M1)]
Each i-root containing $x$ is Moufang (and the corresponding root group is a subgroup of $E$). 
\item[(M2)]
A redundant property, since it follows from (MSTGQ1).
\item[(M3)]
No line $U \I x$ is the unique center of a triad $\{V,W,X\}$ with $V \I x$.
\end{itemize}

He then showed that all flock quadrangles are MSTGQs, and defined his now famous ``Property (G)'' \cite{STGQ89}.

A combination of \cite{flocks87} and \cite{STGQ89} eventually led to the discovery of most of the presently known examples of GQs of order $(q,q^2)$.

Recall that all known finite generalized quadrangles except the Hermitian quadrangles in projective $4$-space have the property that up to a combination of duality and Payne-integration, they are STGQs. This observation follows from the main result of \cite{QKan}.

Besides fundamental work in especially \cite{STGQ89}, not much is known in the classification theory for STGQs.

By results of Hachenberger \cite{Hach} and independently Chen (unpublished) | see \cite{LEGQ}, what we {\em do} know is: 

\begin{theorem}[Chen/Hachenberger \cite{Hach}, see also \cite{LEGQ}]
The parameters of a finite STGQ are always power of the same prime.
\end{theorem}

\bigskip
\section{An observation on flock quadrangles}

\subsection{General Heisenberg groups}

The {\em general Heisenberg group} $\hH_n(q)$ (sometimes also written as $\mathcal{H}_n$ if we do not want to specify $q$) of dimension $2n + 1$ over $\mathbb{F}_q$, with $n$ a positive integer,
is the group of square $(n + 2)\times(n + 2)$-matrices with entries in $\mathbb{F}_q$, of the following form (and with the usual matrix multiplication):

\[ \left(
 \begin{array}{ccc}
 1 & \alpha & c\\
 0 & \mathrm{id}_{n\times n} & \beta^T\\
 0 & 0 & 1\\
 \end{array}
 \right),                                         \]

 \noindent
 where $\alpha, \beta \in \mathbb{F}_q^n$, $c \in \mathbb{F}_q$ and with $\mathrm{id}_{n\times n}$  being the $(n\times n)$-identity matrix.
The group $\hH_n$ is isomorphic to the group $\{(\alpha,c,\beta) \vert \alpha,\beta \in \mathbb{F}_q^n, c \in \mathbb{F}_q\}$, where the group operation $\circ$ is given by
$(\alpha,c,\beta)\circ(\alpha',c',\beta') = (\alpha + \alpha', c + c' + \alpha{\beta'}^T, \beta + \beta')$ (here, $(\cdot)^T$ is a notation for transposition).
The following properties hold for $\hH_n$ (defined over $\mathbb{F}_q$).

\begin{itemize}
\item
$\hH_n$ has exponent $p$ if $q = p^h$ with $p$ an odd prime; it has exponent $4$ if $q$ is even.
\item
The center of $\hH_n$ is given by $Z = Z(\hH_n) = \{(0,c,0) \vert c \in \mathbb{F}_q\}$.         
\item
$[\hH_n,\hH_n] = Z = \Phi(\hH_n)$ and $\hH_n$ is nilpotent of class $2$ ($\Phi(\hH_n)$ is the {\em Frattini subgroup} of $H$, that is, the intersection of all its maximal subgroups).
\end{itemize}

\subsection{Forms and spaces}
\label{fsp}

Let $V$ be the elementary abelian $p$-group $\mathcal{H}_2(q)/Z$.
The map $\chi$ 

\[  \chi: V \times V \mapsto \mathbb{F}_q: (aZ,bZ) \mapsto [a,b]    \]
\noindent
naturally defines a non-singular bilinear alternating form over $\mathbb{F}_q \equiv Z$. So $V$ can be seen as a $4$-dimensional space over $\mathbb{F}_q$, and in the corresponding projective $3$-space over $\mathbb{F}_q$, $\chi$ defines a symplectic polar space $\hW_3(q)$ of rank $2$ (projective index $1$). Here, $\hW_3(q)$ can be defined as the generalized quadrangle which arises as the points
of $\PG(3,q)$, $\mathbb{F}_q$ a finite field,
together with the totally isotropic projective lines of a non-degenerate alternating bilinear form on $\PG(3,q)$. All its points are regular elation points.

\subsection{Special groups and flocks}

The importance of alternating forms in STGQ theory can be read from the following theorem. Its proof uses the connection explained in \S\S \ref{fsp}.

\begin{theorem}[\cite{isomflocks}]
Suppose $H$ is a special $p$-group of order $q^{5}$ for which $Z(H)  =  \Phi(H) =  [H,H]$ is elementary abelian of order $q$.
Suppose $H$ admits  a Kantor family of type $(q^2,q)$, and suppose $\chi$ defines a non-singular bilinear alternating form over $\mathbb{F}_q$. Then $H \cong \mathscr{H}_2(q)$, and 
the corresponding generalized quadrangle $\Gamma$ of order $(q^2,q)$ is a flock quadrangle.
 \end{theorem}

Roughly put (see the citation \cite{isomflocks} for more details), if for an STGQ $(\Gamma,x,E)$, $E$ is isomorphic to a general Heisenberg group,
then $\Gamma$ is a flock quadrangle: 
\begin{equation}
\mbox{STGQ} \ \ + \ \ \mbox{Heisenberg} \ \ \equiv \ \ \mbox{flock}.
\end{equation}

So enough structural knowledge of the group in this case, leads to determination of $\Gamma$. This result/idea is one of the main models
for the classification started in \cite{STGQ}.

\bigskip
\section{Results and corollaries}

We are ready to describe several new results in STGQ theory. Some of them settle long-standing open problems.
Proofs will be published elsewhere (see \cite{STGQ}).

\subsection{STGQs of order $(q,q)$, $q$ odd}

The first result was obtained in 2009, and already announced (a.o.) on the 2010 conference ``Combinatorics 2010'' in Verbania. It was also mentioned (without proof) in the proceedings paper of my talk \cite{KTlocal}.
It can be found in the preprint \cite{STGQ}. I explained several proofs of this result in a lecture at the ``Buildings 2012'' conference in M\"{u}nster.

\begin{theorem}
Let $(\Gamma,x,E)$ be an STGQ of order $(q,q)$, $q$ odd. Then $\Gamma$ is isomorphic to $\hW_3(q)$, $x$ is arbitrary and $E$ is 
isomorphic to $\hH_1(q)$. 
\end{theorem}

\subsection{Payne's conjecture}

The next result completely settles Payne's 2004 question mentioned in the first section.

\begin{theorem}
Let $(\Gamma,x,E)$ be a finite STGQ of order $(s,t)$ with distinct elation groups. Then $s = t^2$, $t$ is a power of $2$, and 
$\Gamma \cong \mathcal{H}(3,t^2)$.  (Moreover, $E$ is known.)
\end{theorem}

\subsection{Fix point theory}

It appears that the following property is crucial (``centrality'').

\quad (C): {\em The group of symmetries with center $x$ is a subgroup of the center of $E$}.

If (C) is true for $(\Gamma,x,E)$, a slightly more general version of the Moufang property holds for any i-root containing $x$. This allows one 
to control the situation to quite a far extent \cite{STGQ}; $E$ then comes in a class of abstract groups which share several distinguished properties with general Heisenberg groups.

\begin{observation}
All known finite STGQs have $(C)$.
\end{observation}

In the study of STGQs (and the known examples), a second property arises naturally:

\quad Property (*): {\em Let $(\Gamma,x,E)$ be an STGQ. Let $Y \I x$. Then $\Gamma$ has $(*)$ {\em at $Y$}, if for some $y \I Y$, $y \ne x$, 
$E_y$ is a normal subgroup of $E$. In that case, $E_y$ is independent of the choice of $y$. The STGQ has $(*)$ if it has $(*)$ at every line on $x$.}\\

Payne's MSTGQs always satisfy (*).

In the case where the parameters of a GQ are of type $(s = t^2,t)$, the next unexpected theorem reveals the intimate connection between  (*) and STGQs. 

\begin{theorem}
A finite EGQ of order $(t^2,t)$ with $(*)$ is an STGQ.
\end{theorem}

A careful and elaborate analysis of fixed point structures in EGQs and STGQs, leads to the following theorem.

\begin{theorem}
If a finite STGQ has $(*)$, and its order is not $(t,t)$ if $t$ is even, property $(C)$ is satisfied.
\end{theorem}

An important corollary is that all i-roots on $x$ are Moufang, if (*) is satisfied (and its order is not $(t,t)$ if $t$ is even). This enables one to understand the elation groups much better | see \cite{STGQ}. When the order is $(t,t)$ and $t$ is even, the situation is (very) different, and an altogether different 
approach is needed to attack this case.

\subsection{Generic STGQs}

Unfortunately (depending on the viewpoint), not all STGQs have (*) at every line through the elation point; a first class of 
counter examples (related to Suzuki groups) is displayed in the preprint \cite{STGQ} | see also \S 8.4 of \cite{LEGQ}. In \cite{STGQ} it is conjectured that for a finite STGQ
$(\Gamma,x,E)$, either zero, one or all lines incident with $x$ have (*). 

In \cite{STGQ}, a very general class of STGQs is introduced which do not satisfy (*). 
Let us call $(\Gamma,x,E)$ {\em generic} if (*) is not satisfied.

 Let $(\Gamma,x,E)$ be an STGQ, and let $\Phi := \Phi(E)$ be the Frattini subgroup of $E$. Define a point-line geometry $\Gamma(\Phi)$ as follows. 
Its lines are the $\Phi$-orbits on the lines incident with $x$ (where the trivial orbit $\{x\}$ is excluded); its points are the $\Phi$-orbits
in the set of points not collinear with $x$. A point $u$ is incident with a line $V$ if at least one $\Gamma$-point of the orbit $V$ is collinear
with some point of the orbit $u$. If $u$ then is incident with $V$, it is easy to see that $V$ is surjectively projected on 
$u$, so that for each $L \I x$, a $\Gamma(\Phi)$-point is incident with precisely one line which is a $\Phi$-orbit on $L$. So each $\Gamma(\Phi)$-point
is incident with precisely $t + 1$ $\Gamma(\Phi)$-lines.

Let $(\Gamma,x,E)$ be a generic STGQ of order $(s,t)$, with Kantor family $(\mE,\mE^*)$.  Let $T = E/\mathbb{S}$, and suppose that 
the only extension of $T$ in $E$ is $E$ itself. 
Suppose also that at least one of the following properties is satisfied.
\begin{itemize}
\item[(a)]
$\Gamma(\Phi)$ is a dual partial linear space. 
\item[(b)]
For each $A \ne B \in \mE$, we have that if $K$ is a maximal subgroup of $H$ which does not contain $A$, then 
$\langle A \cap K, B \rangle \ne H$.

%\item[]
\end{itemize}

It is then shown that $(\Gamma,x,E)$ always contains ideal sub STGQs (that is, it contains sub STGQs containing $x$ and all the lines on $x$).
Whence the parameters are always of type $(t^2,t)$ for this type of generic STGQs.

Finally, it is conjectured that an STGQ either has (*), or is generic with the additional properties assumed above.

\end{document}